\documentclass{amsart}
\usepackage{amstext}
\usepackage{amsthm}
\usepackage{amssymb}

\numberwithin{equation}{section}
\numberwithin{figure}{section}
\theoremstyle{plain}
\newtheorem{thm}{Theorem}

\allowdisplaybreaks
\usepackage[sort&compress,square,numbers]{natbib}
\usepackage{tikz}
\usetikzlibrary{matrix,positioning}

\begin{document}

\title[ Some identities of  Eulerian polynomials]{
Some identities of Eulerian polynomials arising from nonlinear differential equations }

\author{Taekyun Kim}
\address{Department of Mathematics, Kwangwoon University, Seoul 139-701, Republic of Korea}
\email{tkkim@kw.ac.kr}

\author{Dae San Kim}
\address{Department of Mathematics, Sogang University, Seoul 121-742, Republic of Korea}
\email{dskim@sogang.ac.kr}

\begin{abstract}
In this paper, we study nonlinear differential equations arising from
Eulerian polynomials and their applications. From our study of nonlinear
differential equations, we derive some new and explicit identities involving Eulerian
and higher-order Eulerian polynomials.
\end{abstract}



\maketitle
\global\long\def\relphantom#1{\mathrel{\phantom{{#1}}}}

\section{Introduction}

The Eulerian polynomials were introduced by L. Euler in his \emph{ Remarques
sur un beau rapport entre les s\'eries des puissances tant directes
que r\'eciproques} in 1749 (first printed in 1765) where he describes
a method of computing values of the zeta function at negative integers
by a precursor of Abel's theorem applied to a divergent series (see
\cite{key-3,key-4,key-5,key-19,key-22}).

As is well known, the Eulerian polynomials, $A_{n}\left(t\right)$,
$\left(n\ge0\right)$, are defined by the generating function
\begin{equation}
\frac{1-t}{e^{x\left(t-1\right)}-t}=e^{A\left(t\right)x}=\sum_{n=0}^{\infty}A_{n}\left(t\right)\frac{x^{n}}{n!},\quad\left(\text{see \cite{key-14}}\right),\label{eq:1}
\end{equation}
with the usual convention about replacing $A^{n}\left(t\right)$ by
$A_{n}\left(t\right)$.

From (\ref{eq:1}), we can derive the following recurrence relation
for the Eulerian polynomials:
\begin{equation}
\left(A\left(t\right)+\left(t-1\right)\right)^{n}-tA_{n}\left(t\right)=\left(1-t\right)\delta_{0,n},\quad\left(n\ge0\right),\quad\left(\text{see \cite{key-14}}\right).\label{eq:2}
\end{equation}

By (\ref{eq:2}), we easily get
\begin{equation}
A_{0}\left(t\right)=1,\quad A_{n}\left(t\right)=\frac{1}{t-1}\sum_{l=0}^{n-1}\binom{n}{l}A_{l}\left(t\right)\left(t-1\right)^{n-l},\quad\left(n\ge1\right).\label{eq:3}
\end{equation}

Furthermore,
\begin{equation}
\frac{A_{n}\left(t\right)}{\left(1-t\right)^{n+1}}=\sum_{j=0}^{\infty}t^{j}\left(j+1\right)^{n},\quad\left(n\ge0\right),\quad\left(\text{see \cite{key-2,key-3,key-4,key-5,key-6}}\right)\label{eq:4}
\end{equation}

The first few Eulerian polynomials are
\begin{align*}
1+t+t^{2}+t^{3}+\cdots & =\frac{1}{1-t}=\frac{A_{0}\left(t\right)}{1-t},\\
1+2t+3t^{2}+4t^{3}+\cdots & =\frac{1}{\left(1-t\right)^{2}}=\frac{A_{1}\left(t\right)}{\left(1-t\right)^{2}},\\
1+2^{2}t+3^{2}t^{2}+4^{2}t^{3}+\cdots & =\frac{1+t}{\left(1-t\right)^{3}}=\frac{A_{2}\left(t\right)}{\left(1-t\right)^{3}}.
\end{align*}

Recently, several authors has studied some interesting extensions
and modifications of Eulerian polynomials along with related combinatorial,
probabilistic and statistical applications (see \cite{key-1,key-2,key-3,key-4,key-5,key-6,key-7,key-8,key-9,key-10,key-11,key-12,key-13,key-14,key-15,key-16,key-17,key-18,key-19,key-20,key-21,key-22}).

In \cite{key-15}, Kim has studied nonlinear differential equations
arising from Frobenius-Euler numbers and polynomials.

In this paper, we give some new and explicit identities on Eulerian and higher-order Eulerian
polynomials which are derived from solutions of nonlinear differential
equations.

\section{Nonlinear differential equations arising from Eulerian polynomials}

Let us put
\begin{equation}
F=F\left(t,x\right)=\frac{1}{e^{x\left(t-1\right)}-t},\quad\left(t\neq1\right).\label{eq:5}
\end{equation}

Now, we consider the differentiation of $F$ with respect to $x$
while $t$ is being fixed.

\begin{align}
F^{\left(1\right)} & =\frac{d}{dx}F\left(t,x\right)\label{eq:6}\\
 & =\frac{\left(-1\right)e^{x\left(t-1\right)}}{\left(e^{x\left(t-1\right)}-t\right)^{2}}\left(t-1\right)\nonumber \\
 & =\left(1-t\right)\frac{1}{\left(e^{x\left(t-1\right)}-t\right)^{2}}\left(e^{x\left(t-1\right)}-t+t\right)\nonumber \\
 & =\left(1-t\right)\left(F+tF^{2}\right).\nonumber
\end{align}

Thus, by (\ref{eq:6}), we easily get
\begin{align}
F^{\left(2\right)} & =\frac{d}{dx}F^{\left(1\right)}\label{eq:7}\\
 & =\left(1-t\right)\left(F^{\left(1\right)}+2tFF^{\left(1\right)}\right)\nonumber \\
 & =\left(1-t\right)\left(1+2tF\right)F^{\left(1\right)}\nonumber \\
 & =\left(1-t\right)^{2}\left(1+2tF\right)\left(F+tF^{2}\right)\nonumber \\
 & =\left(1-t\right)^{2}\left(F+3tF^{2}+2t^{2}F^{3}\right),\nonumber
\end{align}
and
\begin{align}
F^{\left(3\right)} & =\frac{d}{dx}F^{\left(2\right)}\label{eq:8}\\
 & =\left(1-t\right)^{2}\left(F^{\left(1\right)}+6tFF^{\left(1\right)}+6t^{2}F^{2}F^{\left(1\right)}\right)\nonumber \\
 & =\left(1-t\right)^{2}\left(1+6tF+6t^{2}F^{2}\right)F^{\left(1\right)}\nonumber \\
 & =\left(1-t\right)^{3}\left(1+6tF+6t^{2}F^{2}\right)\left(F+tF^{2}\right)\nonumber \\
 & =\left(1-t\right)^{3}\left(F+7tF^{2}+12t^{2}F^{3}+6t^{3}F^{4}\right).\nonumber
\end{align}

Continuing this process, we set
\begin{align}
F^{\left(N\right)} & =\left(\frac{d}{dx}\right)^{N}F\left(t,x\right)\label{eq:9}\\
 & =\left(\frac{d}{dx}\right)^{N}\left(\frac{1}{e^{x\left(t-1\right)}-t}\right)\nonumber \\
 & =\left(1-t\right)^{N}\sum_{i=1}^{N+1}a_{i-1}\left(N,t\right)F^{i},\quad\left(N\in\mathbb{N}\cup\left\{ 0\right\}\right).\nonumber
\end{align}

From \ref{eq:9}, we can derive the following equation (\ref{eq:10}):
\begin{align}
F^{\left(N+1\right)} & =\frac{d}{dx}F^{\left(N\right)}\label{eq:10}\\
 & =\left(1-t\right)^{N}\sum_{i=1}^{N+1}a_{i-1}\left(N,t\right)iF^{i-1}F^{\left(1\right)}\nonumber \\
 & =\left(1-t\right)^{N+1}\sum_{i=1}^{N+1}a_{i-1}\left(N,t\right)iF^{i-1}\left(F+tF^{2}\right)\nonumber \\
 & =\left(1-t\right)^{N+1}\left\{ \sum_{i=1}^{N+1}a_{i-1}\left(N,t\right)iF^{i}+\sum_{i=1}^{N+1}a_{i-1}\left(N,t\right)itF^{i+1}\right\} \nonumber \\
 & =\left(1-t\right)^{N+1}\left\{ \sum_{i=1}^{N+1}a_{i-1}\left(N,t\right)iF^{i}+\sum_{i=2}^{N+2}a_{i-2}\left(N,t\right)\left(i-1\right)tF^{i}\right\} \nonumber \\
 & =\left(1-t\right)^{N+1}\left\{ a_{0}\left(N,t\right)F+\left(N+1\right)t a_{N}\left(N,t\right)F^{N+2}\right.\nonumber \\
 & \relphantom =\left.+\sum_{i=2}^{N+1}\left(ia_{i-1}\left(N,t\right)+\left(i-1\right)t a_{i-2}\left(N,t\right)\right)F^{i}\right\}
 .\nonumber
\end{align}

By replacing $N$ by $N+1$ in (\ref{eq:9}), we get
\begin{equation}
F^{\left(N+1\right)}=\left(1-t\right)^{N+1}\sum_{i=1}^{N+2}a_{i-1}\left(N+1,t\right)F^{i}.\label{eq:11}
\end{equation}

From (\ref{eq:10}) and (\ref{eq:11}), we can derive the following
recurrence relation for the coefficients $a_{i}\left(N,t\right)$:
\begin{align}
a_{0}\left(N+1,t\right) & =a_{0}\left(N,t\right),\label{eq:12}\\
a_{N+1}\left(N+1,t\right) & =\left(N+1\right)ta_{N}\left(N,t\right),\label{eq:13}
\end{align}
and
\begin{equation}
a_{i-1}\left(N+1,t\right)=\left(i-1\right)ta_{i-2}\left(N,t\right)+ia_{i-1}\left(N,t\right),\label{eq:14}
\end{equation}
where $2\le i\le N+1$.

It is not difficult to show that
\begin{equation}
F=F^{\left(0\right)}=a_{0}\left(0,t\right)F.\label{eq:15}
\end{equation}

Thus, by (\ref{eq:15}), we have
\begin{equation}
a_{0}\left(0,t\right)=1.\label{eq:16}
\end{equation}

From (\ref{eq:7}) and (\ref{eq:9}), we note that
\begin{align}
\left(1-t\right)\left(F+tF^{2}\right) & =F^{\left(1\right)}\label{eq:17}\\
 & =\left(1-t\right)\sum_{i=1}^{2}a_{i-1}\left(1,t\right)F^{i}\nonumber \\
 & =\left(1-t\right)\left\{ a_{0}\left(1,t\right)F+a_{1}\left(1,t\right)F^{2}\right\} .\nonumber
\end{align}

By comparing the coefficients on both sides of (\ref{eq:17}),
we have
\begin{equation}
a_{0}\left(1,t\right)=1,\quad a_{1}\left(1,t\right)=t.\label{eq:18}
\end{equation}

From (\ref{eq:18}), we note that
\begin{equation}
a_{0}\left(N+1,t\right)=a_{0}\left(N,t\right)=\cdots=a_{0}\left(1,t\right)=a_{0}\left(0,t\right)=1,\label{eq:19}
\end{equation}
and
\begin{align}
a_{N+1}\left(N+1,t\right) & =\left(N+1\right)ta_{N}\left(N,t\right)\label{eq:20}\\
 & =\left(N+1\right)tNta_{N-1}\left(N-1,t\right)\nonumber \\
 & =t^{2}\left(N+1\right)Na_{N-1}\left(N-1,t\right)\nonumber \\
 & \vdots\nonumber \\
 & =t^{N}\left(N+1\right)N\cdots2a_{1}\left(1,t\right)\nonumber \\
 & =t^{N+1}\left(N+1\right)!.\nonumber
\end{align}

So, we have the matrix $\left(a_{i}\left(j,t\right)\right)_{0\le i,j\le N}$
as follows:

\begin{equation*}
 \begin{tikzpicture}[baseline=(current  bounding  box.west)]
  \matrix (mymatrix) [matrix of math nodes,left delimiter={[},right
delimiter={]}]
  {
    1  & 1   & 1 & 1 & \cdots & 1 \\
      & 1!t   &  &  &    \\
     &  &  2!t^2  &  &     \\
     &  & & 3!t^3 & &  \\
    &  &  & &\ddots &  \\
    &  &  & & & N!t^N \\
  };
\node[xshift=-57pt,yshift=55pt] {$0$};
\node[xshift=-42pt,yshift=55pt] {$1$};
\node[xshift=-22pt,yshift=55pt] {$2$};
\node[xshift=1pt,yshift=55pt] {$3$};
\node[xshift=45pt,yshift=55pt] {$N$};
\node[xshift=-78pt,yshift=40pt] {$0$};
\node[xshift=-78pt,yshift=26pt] {$1$};
\node[xshift=-78pt,yshift=11pt] {$2$};
\node[xshift=-78pt,yshift=-4pt] {$3$};
\node[xshift=-78pt,yshift=-40pt] {$N$};
\node[xshift=-40pt,yshift=-30pt] {\LARGE$0$};
\end{tikzpicture}\label{eq:27}
\end{equation*}

From (\ref{eq:14}), we have
\begin{align}
a_{1}\left(N+1,t\right) & =ta_{0}\left(N,t\right)+2a_{1}\left(N,t\right)\label{eq:21}\\
 & =ta_{0}\left(N,t\right)+2\left\{ ta_{0}\left(N-1,t\right)+2a_{1}\left(N-1,t\right)\right\} \nonumber \\
 & =t\left\{ a_{0}\left(N,t\right)+2a_{0}\left(N-1,t\right)\right\} +2^{2}\left\{ ta_{0}\left(N-2,t\right)+2a_{1}\left(N-2,t\right)\right\} \nonumber \\
 & =t\left\{ a_{0}\left(N,t\right)+2a_{0}\left(N-1,t\right)+2^{2}a_{0}\left(N-2,t\right)\right\} +2^{3}a_{1}\left(N-2,t\right)\nonumber \\
 & \vdots\nonumber \\
 & =t\sum_{i=0}^{N-1}2^{i}a_{0}\left(N-i,t\right)+2^{N}a_{1}\left(1,t\right)\nonumber \\
 & =t\sum_{i=0}^{N}2^{i}a_{0}\left(N-i,t\right),\nonumber
\end{align}
\begin{align}
a_{2}\left(N+1,t\right) & =2ta_{1}\left(N,t\right)+3a_{2}\left(N,t\right)\label{eq:22}\\
 & =2ta_{1}\left(N,t\right)+3\left\{ 2ta_{1}\left(N-1,t\right)+3a_{2}\left(N-1,t\right)\right\} \nonumber \\
 & =2t\left\{ a_{1}\left(N,t\right)+3a_{1}\left(N-1,t\right)\right\} +3^{2}\left\{ 2ta_{1}\left(N-2,t\right)+3a_{2}\left(N-2,t\right)\right\} \nonumber \\
 & =2t\left\{ a_{1}\left(N,t\right)+3a_{1}\left(N-1,t\right)+3^{2}a_{1}\left(N-2,t\right)\right\} +3^{3}a_{2}\left(N-2,t\right)\nonumber \\
 & \vdots\nonumber \\
 & =2t\sum_{i=0}^{N-2}3^{i}a_{1}\left(N-i,t\right)+3^{N-1}a_{2}\left(2,t\right)\nonumber \\
 & =2t\sum_{i=0}^{N-1}3^{i}a_{1}\left(N-i,t\right),\nonumber
\end{align}
and
\begin{align}
a_{3}\left(N+1,t\right) & =3ta_{2}\left(N,t\right)+4a_{3}\left(N,t\right)\label{eq:23}\\
 & =3ta_{2}\left(N,t\right)+4\left\{ 3ta_{2}\left(N-1,t\right)+4a_{3}\left(N-1,t\right)\right\} \nonumber \\
 & =3t\left\{ a_{2}\left(N,t\right)+4a_{2}\left(N-1,t\right)\right\} +4^{2}\left\{ 3ta_{2}\left(N-2,t\right)+4a_{3}\left(N-2,t\right)\right\} \nonumber \\
 & =3t\left\{ a_{2}\left(N,t\right)+4a_{2}\left(N-1,t\right)+4^{2}a_{2}\left(N-2,t\right)\right\} +4^{3}a_{3}\left(N-2,t\right)\nonumber \\
 & \vdots\nonumber \\
 & =3t\sum_{i=0}^{N-3}4^{i}a_{2}\left(N-i,t\right)+4^{N-2}a_{3}\left(3,t\right)\nonumber \\
 & =3t\sum_{i=0}^{N-2}4^{i}a_{2}\left(N-i,t\right).\nonumber
\end{align}

Continuing this process, we get
\begin{equation}
a_{j}\left(N+1,t\right)=jt\sum_{i=0}^{N-j+1}\left(j+1\right)^{i}a_{j-1}\left(N-i,t\right),\quad\left(1\le j\le N+1\right).\label{eq:24}
\end{equation}

Therefore, by (\ref{eq:24}), we obtain the following theorem.
\begin{thm}
\label{thm:1} For each fixed $t\left(\neq1\right)$ and $N\in\mathbb{N}\cup\left\{ 0\right\}$,
$F=F\left(t,x\right)=\frac{1}{e^{x\left(t-1\right)}-t}$ satisfies
the nonlinear differential equation
\begin{equation}
\left(\frac{d}{dx}\right)^{N}F=\left(1-t\right)^{N}\sum_{i=1}^{N+1}a_{i-1}\left(N,t\right)F^{i},\label{eq:25}
\end{equation}
where $a_{0}\left(N,t\right)=a_{0}\left(N-1,t\right)=\cdots=a_{0}\left(1,t\right)=a_{0}\left(0,t\right)=1$,
\[
a_{i}\left(N,t\right)=it\sum_{j=0}^{N-i}\left(i+1\right)^{j}a_{i-1}\left(N-j-1,t\right)\quad\left(1\le j\le N\right).
\]
 \end{thm}

Taking the $N$-th derivative with respect to $x$ on both sides
of (\ref{eq:1}), we obtain
\begin{align}
\left(\frac{d}{dx}\right)^{N}\frac{1-t}{e^{x\left(t-1\right)}-t} & =\sum_{n=N}^{\infty}A_{n}\left(t\right)\left(n\right)_{N}\frac{x^{n-N}}{n!}\label{eq:28}\\
 & =\sum_{n=0}^{\infty}A_{n+N}\left(t\right)\left(n+N\right)_{N}\frac{x^{n}}{\left(n+N\right)!}\nonumber \\
 & =\sum_{n=0}^{\infty}A_{n+N}\left(t\right)\frac{x^{n}}{n!}.\nonumber
\end{align}

On the other hand, from (\ref{eq:25}), we have
\begin{align}
 & \left(\frac{d}{dx}\right)^{N}\left(\frac{1-t}{e^{x\left(t-1\right)}-1}\right)\label{eq:29}\\
 & =\left(1-t\right)^{N+1}\sum_{i=1}^{N+1}a_{i-1}\left(N,t\right)\left(1-t\right)^{-i}\left(\frac{1-t}{e^{x\left(t-1\right)}-t}\right)^{i}\nonumber \\
 & =\sum_{i=1}^{N+1}a_{i-1}\left(N,t\right)\left(1-t\right)^{N+1-i}\sum_{n=0}^{\infty}A_{n}^{\left(i\right)}\left(t\right)\frac{x^{n}}{n!}\nonumber \\
 & =\sum_{n=0}^{\infty}\left(\sum_{i=1}^{N+1}a_{i-1}\left(N,t\right)\left(1-t\right)^{N+1-i}A_{n}^{\left(i\right)}\left(t\right)\right)\frac{x^{n}}{n!},\nonumber
\end{align}
where $A_{n}^{\left(i\right)}\left(t\right)$ are called the higher-order
Eulerian polynomials and defined by the generating function
\begin{equation}
\left(\frac{1-t}{e^{x\left(t-1\right)}-1}\right)^{m}=\sum_{n=0}^{\infty}A_{n}^{\left(m\right)}\left(t\right)\frac{x^{n}}{n!}.\label{eq:30}
\end{equation}

Therefore, by (\ref{eq:28}) and (\ref{eq:29}), we obtain the following
theorem.
\begin{thm}
\label{thm:2} For all $t$, and $n, N\in\mathbb{N}\cup\left\{ 0\right\}$,
we have
\[
A_{n+N}\left(t\right)=\sum_{i=1}^{N+1}a_{i-1}\left(N,t\right)\left(1-t\right)^{N+1-i}A_{n}^{\left(i\right)}\left(t\right).
\]

\end{thm}
Note here that, as both sides are polynomials and it holds for all $t\left(\neq1\right)$,
it is true as polynomials. Explicit expressions for $a_{i}\left(N,t\right)$,
$\left(1\le i\le N\right)$, are given by
\begin{align}
a_{1}\left(N,t\right) & =t\sum_{j=0}^{N-1}2^{j}a_{0}\left(N-j-1,t\right)\label{eq:31}\\
 & =t\sum_{j=0}^{N-1}2^{j}\nonumber \\
 & =t\left(2^{N}-1\right),\nonumber
\end{align}
\begin{align}
a_{2}\left(N,t\right) & =2t\sum_{j_{1}=0}^{N-2}3^{j_{1}}a_{1}\left(N-j_{1}-1,t\right)\label{eq:32}\\
 & =2t\sum_{j_{1}=0}^{N-2}3^{j_{1}}t\left(2^{N-j_{1}-1}-1\right)\nonumber \\
 & =2t^{2}\sum_{j_{1}=0}^{N-2}3^{j_{1}}\left(2^{N-j_{1}-1}-1\right),\nonumber
\end{align}
\begin{align}
a_{3}\left(N,t\right) & =3t\sum_{j_{2}=0}^{N-3}4^{j_{2}}a_{2}\left(N-j_{2}-1,t\right)\label{eq:33}\\
 & =3t\sum_{j_{2}=0}^{N-3}4^{j_{2}}2t^{2}\sum_{j_{1}=0}^{N-j_{2}-3}3^{j_{1}}\left(2^{N-j_{2}-j_{1}-2}-1\right)\nonumber \\
 & =3!t^{3}\sum_{j_{2}=0}^{N-3}\sum_{j_{1}=0}^{N-j_{2}-3}4^{j_{2}}3^{j_{1}}\left(2^{N-j_{2}-j_{1}-2}-1\right),\nonumber
\end{align}
and
\begin{align}
a_{4}\left(N,t\right) & =4t\sum_{j_{3}=0}^{N-4}5^{j_{3}}a_{3}\left(N-j_{3}-1,t\right)\label{eq:34}\\
 & =4t\sum_{j_{3}=0}^{N-4}5^{j_{3}}3!t^{3}\sum_{j_{2}=0}^{N-j_{3}-4}\sum_{j_{1}=0}^{N-j_{3}-j_{2}-4}4^{j_{2}}3^{j_{1}}\left(2^{N-j_{3}-j_{2}-j_{1}-3}-1\right)\nonumber \\
 & =4!t^{4}\sum_{j_{3}=0}^{N-4}\sum_{j_{2}=0}^{N-j_{3}-4}\sum_{j_{1}=0}^{N-j_{3}-j_{2}-4}5^{j_{3}}4^{j_{2}}3^{j_{1}}\left(2^{N-j_{3}-j_{2}-j_{1}-3}-1\right).\nonumber
\end{align}

Continuing this process, we have
\begin{align}
a_{i}\left(N,t\right) & =i!t^{i}\sum_{j_{i-1}=0}^{N-i}\sum_{j_{i-2}=0}^{N-j_{i-1}-i}\cdots\sum_{j_{1}=0}^{N-j_{i-1}-\cdots-j_{2}-i}\left(i+1\right)^{j_{i-1}}i^{j_{i-2}}\cdots3^{j_{1}}\label{eq:35}\\
 & \relphantom =\times\left(2^{N-j_{i-1}-j_{i-2}-\cdots-j_{1}-i+1}-1\right),\quad\left(1\le i\le N\right).\nonumber
\end{align}

Recall that
\begin{equation}
\frac{A_{n}\left(t\right)}{\left(1-t\right)^{n+1}}=\sum_{j=0}^{\infty}t^{j}\left(j+1\right)^{n},\quad\left(n\ge0\right).\label{eq:36}
\end{equation}

Thus, by (\ref{eq:36}), we get
\begin{align}
 & \sum_{j=0}^{\infty}t^{j}\left(j+1\right)^{n+N}\label{eq:37}\\
 & =\frac{A_{n+N}\left(t\right)}{\left(1-t\right)^{n+N+1}}\nonumber \\
 & =\sum_{i=1}^{N+1}a_{i-1}\left(N,t\right)\left(1-t\right)^{-n-i}A_{n}^{\left(i\right)}\left(t\right)\nonumber \\
 & =\left(1-t\right)^{-n}\sum_{i=1}^{N+1}a_{i-1}\left(N,t\right)\left(1-t\right)^{-i}A_{n}^{\left(i\right)}\left(t\right).\nonumber
\end{align}

Therefore, by (\ref{eq:37}), we obtain the following theorem.
\begin{thm}
\label{thm:3} For $n, N\in\mathbb{N}\cup\left\{ 0\right\} $,
we have
\[
\sum_{j=0}^{\infty}t^{j}\left(j+1\right)^{n+N}=\left(1-t\right)^{-n}\sum_{i=1}^{N+1}a_{i-1}\left(N,t\right)\left(1-t\right)^{-i}A_{n}^{\left(i\right)}\left(t\right).
\]

\end{thm}
\bibliographystyle{amsplain}
\providecommand{\bysame}{\leavevmode\hbox to3em{\hrulefill}\thinspace}
\providecommand{\MR}{\relax\ifhmode\unskip\space\fi MR }
\providecommand{\MRhref}[2]{%
  \href{http://www.ams.org/mathscinet-getitem?mr=#1}{#2}
}
\providecommand{\href}[2]{#2}

\end{document}